\documentclass[11pt,leqno]{article}
\usepackage{amssymb,amsmath,amsthm}
\makeatletter
\def\@fnsymbol#1{\ensuremath{\ifcase#1\or 1\or 2\fi}}
\makeatother

\newcommand{\car}{{\rm char}}
\newcommand{\Int}{{\rm Int}}

\newcommand{\eq}{\begin{eqnarray}}
\newcommand{\eneq}{\end{eqnarray}}
\newcommand{\eqn}{\begin{eqnarray*}}
\newcommand{\eneqn}{\end{eqnarray*}}

\newcommand{\indlim}[1][]{\mathop{\varinjlim}\limits_{#1}}

\def\etens{\mathbin{\boxtimes}}
\newcount\temp
\temp=\catcode`\@
\catcode`\@=11
\def\@bletens{\mathbin{\etens^{L}}}
\def\@letens_#1{\mathbin{\etens_{\raise1.5ex\hbox to-.1em{}#1}^{L}}}
\def\letens{\@ifnextchar _{\@letens}{\@bletens}}
\catcode`\@=\temp

\def\phi{{\varphi}}
\def\epsilon{\varepsilon}

\newcommand{\md[1]}{{\rm Mod}({#1})}

\newcommand{\mdco[1]}{{\rm Mod_{coh}}({#1})}

\newcommand{\II[1]}{{\rm Ind}({#1})}

\def\shb{\mathcal{B}}

\def\shd{\mathcal{D}}
\def\she{\mathcal{E}}

\def\shh{\mathcal{H}}

\def\shj{\mathcal{J}}

\def\shl{\mathcal{L}}
\def\shm{\mathcal{M}}
\def\shn{\mathcal{N}}
\def\sho{\mathcal{O}}

\def\shs{\mathcal{S}}

\newcommand{\C}{\mathbb{C}}

\newcommand{\R}{\mathbb{R}}

\newcommand{\Z}{\mathbb{Z}}

\DeclareMathOperator{\coker}{coker}

\renewcommand{\to}[1][]{\xrightarrow[#1]{}}
\newcommand{\from}[1][]{\xleftarrow[#1]{}}
\newcommand{\isoto}[1][]{\xrightarrow[#1]{\sim}}

\newcommand{\Endo}[1][]{\mathrm{End}_{\raise1.5ex\hbox to.1em{}#1}}

\newcommand{\Hom}[1][]{\mathrm{Hom}_{\raise1.5ex\hbox to.1em{}#1}}

\newcommand{\RHom}[1][]{\mathrm{RHom}_{\raise1.5ex\hbox to.1em{}#1}}

\newcommand{\Ext}[2][]{\mathrm{Ext}_{\raise1.5ex\hbox to.1em{}#1}^{#2}}

\newcommand{\THom}[1][]{\mathrm{THom}_{\raise1.5ex\hbox to.1em{}#1}}

\newcommand{\tens}[1][]{\mathbin{\otimes_{\raise1.5ex\hbox to-.1em{}#1}}}

\newcommand{\ltens}[1][]{\mathbin{\otimes_{\raise1.5ex\hbox to-.1em{}#1}^{L}}}

\newcommand{\reim}[1]{{R#1}_!}

\newcommand{\opb}[1]{#1^{-1}}

\renewcommand{\hom}[1][]{{\mathcal{H}om}_{\raise1.5ex\hbox to.1em{}#1}}

\newcommand{\rhom}[1][]{{R\mathcal{H}om}_{\raise1.5ex\hbox to.1em{}#1}}

\newcommand{\thom}[1][]{{T\mathcal{H}om}_{\raise1.5ex\hbox to.1em{}#1}}
\newcommand{\tmuhom}[1][]{{T\mu\mathcal{H}om}_{\raise1.5ex\hbox to.1em{}#1}}
\newcommand{\plus}{{\widehat{+}}}
\newcommand{\sect}{\Gamma}
\newcommand{\rsect}{\mathrm{R}\Gamma}
\newtheorem{theorem}{Theorem}[section]

\newtheorem{lemma}[theorem]{Lemma}
\newtheorem{corollary}[theorem]{Corollary}

\theoremstyle{definition}

\newtheorem{definition}[theorem]{Definition}

\newenvironment{nnum}{
  \begin{enumerate}
  \itemsep=0pt
  
  }
  {\end{enumerate}}

\newenvironment{anum}{
  \begin{enumerate}
  \itemsep=0pt
  
  }
  {\end{enumerate}}

\newcommand{\bnum}{\begin{nnum}}
\newcommand{\enum}{\end{nnum}}
\newcommand{\banum}{\begin{anum}}
\newcommand{\eanum}{\end{anum}}
\newcommand{\ol}{\overline}

\newcommand{\bysame}{\rule{2.5em}{.1mm}\,}

\usepackage[all]{xy}

\CompileMatrices

\author{
Masaki Kashiwara
\thanks{The first named author benefits from
a ``Chaire Internationale de Recherche Blaise Pascal de l'Etat et de la
R\'egion d'Ile-de-France, g\'er\'ee par la
Fondation de l'Ecole Normale Sup\'erieure''.}\enspace ,
Teresa Monteiro Fernandes
\thanks{Research supported by FCT and Programa Ci\^encia,
Tecnologia e Inova\c c\~ao do Quadro Comunit\'ario de Apoio}\enspace, 
and Pierre Schapira}

\title{Micro-support and Cauchy problem for temperate solutions 
of regular $\cal D$-Modules}
\date{December 2000}

\begin{document}

\maketitle
\begin{abstract}
Let $X$ be a complex manifold, $V$ a smooth involutive
submanifold of $T^*X$, $\cal M$ a microdifferential system regular
along $V$, 
 and  $F$  an
$\R$-constructible sheaf on $X$. We study the complex of temperate
microfunction solutions of $\shm$ associated with $F$, that is, the
complex 
$\rhom[\shd_X](\shm,{\cal {T}}\mu hom(F,\sho_X))$. We give  a bound to
its micro-support and solve the Cauchy problem 
 under a suitable hyperbolicity assumption.
\end{abstract}

\section{Introduction}\label{section:intro}
The Cauchy problem for solutions of linear differential operators 
as well as the problem of propagation of singularities, 
are closely related subjects which have been intensively studied
in the 80'th. In the analytic case, it is shown in \cite{K-S1} that
these problems may be reduced to purely geometric ones, using sheaf
theory, the only analytic tool being the Cauchy-Kowalevski
theorem. 

To be more precise, recall that a system of linear partial
differential operators on a complex
manifold $X$ is the data of a coherent module $\shm$ over the sheaf of
rings $\shd_X$ of holomorphic differential operators. 
Let $F$ be a complex of sheaves on $X$ with $\R$-constructible cohomologies
(one says an $\R$-constructible sheaf, for short).
The complex of ``generalized functions'' associated with $F$ is
described by the complex $\rhom(F,\sho_X)$, 
and the complex of solutions of $\shm$
with values in this complex is described by the complex
$$
\rhom[\shd_X](\shm,\rhom(F,\sho_X)).
$$

One may also microlocalize the problem by replacing
$\rhom(F,\sho_X)$ with $\mu hom(F,\sho_X)$. In (loc cit) one shows that
most of the properties of this complex, especially those related to
propagation or Cauchy problem,  are encoded in two geometric objects,
both living in the cotangent bundle $T^*X$, the characteristic variety
of the system
$\shm$, denoted by $\car(\shm)$, and the micro-support of $F$, denoted
by $SS(F)$.

The complex $\rhom(F,\sho_X)$ allows us to treat various situations.
For example if $M$ is a real manifold and $X$ is a complexification of $M$,
by taking as $F$ the dual $D'(\C_M)$ of the constant sheaf on $M$, 
one obtains the sheaf $\shb_M$ of Sato's hyperfunctions. If $Z$ is a complex
analytic hypersurface of $X$ and $F=\C_Z[-1]$ is the (shifted)
constant sheaf on $Z$, one obtains 
the sheaf of holomorphic functions with singularities on $Z$.

\hyphenation{Sch-wartz}

However, the complex $\rhom(F,\sho_X)$ does not allow us to 
treat sheaves associated with holomorphic functions  with
temperate growth such as Schwartz's distributions or meromorphic
functions with poles on $Z$. To consider such cases, one has to
replace it by the complex  ${\cal T}hom(F,\sho_X)$
of temperate cohomology, introduced in \cite{K1} or its
microlocalization ${\cal T}\mu hom(F, {\cal O}_X)$ constructed by
Andronikof \cite {A}. At this stage, a serious difficulty appears: the
geometric methods of
\cite{K-S1} do not apply any more, and indeed, it is a well known fact
that to solve for example the
Cauchy problem for distributions requires more informations 
than the data of the characteristic
variety of the system.

In fact, very little is known
concerning the problems of propagation of singularities and
the Cauchy problem in the space of distributions, apart the 
case where
$\shm$ has real simple characteristics (see \cite{D'A-T,Ko-T}
for a formulation in the language of sheaves and ${\cal D}$-modules)
and some very specific
situation (e.g., operators with constant coefficients on $\R^n$). We
refer to \cite{Ho} for historical and bibliographical comments.

In this paper we give an estimate for the microsupport of the sheaf of 
temperate microfunction solutions associated with an 
$\mathbb{R}$-constructible object $F$, when ${\cal M}$ has regular
singularities along an involutive manifold $V$
in the sense of \cite{K-O}. 
More precisely, we prove the estimate 
\eq\label{eq:1}
&&SS(\rhom[\shd_X](\shm,{\cal T}
\mu hom(F,{\cal O}_X)))\subset \rho_V^{-1}(C_V (SS(F))^a),
\eneq
where $\rho_V\colon V\times_{T^*X} T^*(T^*X)\isoto
V\times_{T^*X} T(T^*X)\to T_V(T^*X)$ is the projection
and
$a$ is the antipodal map of $T_V(T^*X)$
as a vector bundle over $V$.

One can translate this result as follows.
For a bicharacteristic leaf $\Sigma$ of $V$,
one has $\Sigma\times_{T^*X}T_V(T^*X)\simeq T^*\Sigma$,
and $\Sigma\times_{V}C_V (SS(F))^a$ may be regarded as a subset of
$T^*\Sigma$.
Then for any $\Sigma$,
\eqref{eq:1} implies
\eq
&&SS(\rhom[\shd_X](\shm,{\cal T}
\mu hom(F,{\cal O}_X))\vert_\Sigma)\subset \Sigma\times_{V}C_V (SS(F))^a.
\eneq

 What makes this
paper original, in our opinion, is that we treat general
$\R$-constructible sheaves $F$. Let us illustrate our results by an
example. 

 We consider  a smooth morphism $f\colon X\to Y$,  we set
$V=X\times_YT^*Y$, and we assume that
$\shm$ is a coherent module 
regular along $V$.
Let $M$ be a real analytic manifold with complexification
$X$, $S$ a closed subanalytic subset of $M$. We obtain the estimate
$$
SS(\rhom[\shd_X](\shm,\sect_S(\shd b_M)))\subset V\plus (SS(\C_S))^a,
$$
where $\shd b_M$ denote the sheaf of distributions on $M$, 
and the operation $\plus$ 
is defined in \cite{K-S1} and recalled in Section \ref{S:1}.

\section{Notations and main results}\label{S:1}
We will mainly follow the notations in \cite{K-S1}.

Let $X$ be a real analytic manifold. We shall denote by 
$\tau\colon TX\to X$ the tangent bundle to $X$ and
by $\pi\colon T^*X\to X$ the
 cotangent bundle. Set $\dot{T}^*X=T^*X\setminus X$ and   
$\dot{\pi}\colon
\dot{T}^*X\to X$  the projection $\dot{T}^*X\to X$.
For a smooth submanifold $Y$ of $X$, $T_YX$ denotes
 the normal bundle to $Y$ and $T^*_YX$ the conormal bundle. 
In particular, $T^*_X X$ is identified with $X$, the zero section. 

For a submanifold $Y$ of $X$ and a subset  $S$ of $X$, 
we denote by $C_Y(S)$ the normal cone to $S$
along $Y$,  a conic subset of $T_YX$.

If $A$ and $B$ are two conic
subsets of $T^*X$, the operation $A\plus B$ is defined in (loc cit) and
will be recalled below. The set $A^a$ denotes the image of $A$ by the
antipodal map, $(x;\xi)\mapsto (x;-\xi)$. 

For a cone $\gamma\subset TX$,
the polar cone $\gamma^\circ$ to $\gamma$ is the closed convex cone of $T^*X$ 
defined by
$$\gamma^\circ=\{(x;\xi);x\in\pi(\gamma), \langle v,\xi\rangle\geq 0
\,\, \text{for\,\, any}\,\,
(x;v)\in \gamma\}.$$

Let $f\colon X\to Y$ be a morphism of complex manifolds. 
One has two natural morphisms

$$ T^*X\from[f_d]X\times_Y T^*Y\to[f_\pi] T^*Y$$
(In \cite{K-O}, $f_d$ is denoted by ${}^{t}f'$.)

We denote by $D(\C_X)$ (respectively $D^b(\C_X)$,
$D^b_{\R-c}(\C_X)$) the derived category 
of complexes of sheaves of $\C$-vector spaces (respectively 
with bounded cohomologies,  with 
bounded and $\R$-constructible cohomologies).

We denote by $D'_X$ the duality functor on $D^b(\C_X)$, defined by 
$$ D'_X(F)=\rhom(F,\C_X).$$

If $F$ is an object of $D^b(\C_X)$, 
 $SS(F)$ denotes its micro-support, a closed
$\R^{+}$-conic involutive subset of $T^*X$.

On a complex manifold $X$ we consider the sheaf  $\sho_X$  of
holomorphic functions, the sheaf $\Theta_X$ of holomorphic vector
fields, the sheaf $\shd_X$ 
of linear holomorphic differential operators of finite order, 
and its subsheaves $\shd_X(m)$ of operators of order at most $m$. 
We shall also consider the sheaf $\she_X$ on $T^*X$ of 
microdifferential operators of finite order 
(\cite{S-K-K} and \cite{S} for an exposition) 
and its subsheaves $\she_X(m)$ of operators of order at most $m$. 
We denote by 
$\mdco[\shd_X]$ (respectively by 
$\mdco[\she_X]$) the abelian category of  coherent $\shd_X$-modules 
(respectively coherent $\she_X$-modules).
We denote by $D^b(\shd_X)$ the bounded derived category of left
$\shd_X$-modules and by $D^b_{\rm coh}(\shd_X)$ its full triangulated
category consisting of objects with coherent cohomologies. We define
similarly $D^b(\she_X)$ and  $D^b_{\rm coh}(\she_X)$.

The notion of regularity of an $\she_X$-module along an involutive 
submanifold $V$ of $T^*X$ will be recalled in Section \ref{section:reg}.

The main purpose of this paper is to prove the three following results.

Let $V$ be an involutive vector subbundle of $T^*X$.
Let $\shm$ be a coherent $\shd_X$-module regular along $V$
and let $F\in D^b_{\R-c}(\C_X)$.

\begin{theorem}\label{T:1} 
We have the estimate:
\begin{equation}\label{E:2}
SS\big(\rhom[\shd_X](\shm, {\cal T} hom (F,\sho_X))\big)\subset
V{\plus}SS(F)^a. 
\end{equation}
\end{theorem}

Note that $V{\plus}SS(F)^a$ coincides with the
closure of $V{+}SS(F)^a$ by Lemma \ref{L:15}.

As a particular case of Theorem \ref{T:1}, 
assume that X is the complexification
of a real analytic manifold $M$ and let
$S$ be a closed subanalytic subset of $M$. We obtain the estimate:
$$ SS\rhom[\shd_X](\shm,\sect_S(\shd b_M))
\subset (SS(\C_S))^a{\plus} V.$$

For a complex submanifold $Z$ of $X$, we denote by $\shm_Z$ the induced
system of $\shm$ on $Z$.

\begin{theorem}\label{T:2} 
Assume 
$T_Z^*X\cap((SS F)^a {\plus}V)\subset T_X ^*X$.
Then
$$\rhom[\shd_X](\shm, {\cal T} hom(F,\sho_X))\vert_Z 
       \simeq\rhom[\shd_Z](\shm_Z,{\cal T} hom (F_Z, \sho_Z)).$$
\end {theorem}

 Let $U$ be a conic open subset of $\dot{T}^*X$, and  let $V$ be a
closed smooth conic regular involutive submanifold in $U$. We denote by
$$\rho_V : V\times_{T^*X} T^*(T^*X)\isoto
V\times_{T^*X} T(T^*X)\to T_V T^*X$$ the canonical
projection.

 For an  $\mathbb {R}$-constructible sheaf  $F$,
the cohomology of ${\cal T} \mu hom(F,{\cal O}_X)$ is provided with an action of
${\cal{E}}_X$ as proved in \cite{A}, therefore, when ${\cal{T}}\mu hom(F,{\cal
O}_X)$ is concentrated in a single degree, we regard it as an 
${\cal{E}}_X$-module.

\begin{theorem} \label{T:40}
Let $U$ and $V$ be as above and let ${\cal M}$ be a coherent
 ${\cal E}_X$-module defined on
$U$ and regular along $V$. Assume that $F\in D^b_{\R-c}(\C_X)$
and ${\cal{T}}\mu hom(F,{\cal O}_{X})\vert_U$ 
is concentrated in a single degree. Then
\begin{equation}\label{E:52}
 SS(R{\cal H}om_{{\cal E}_X}({\cal M},{\cal{T}}\mu hom(F,{\cal
O}_{X})))\subset {{\rho_{V}}^{-1}(C_{V}(SS(F))^a)}.    
\end{equation}
\end{theorem}

For a bicharacteristic leaf $\Sigma$ of $V$,
one has $\Sigma\times_{T^*X}T_V(T^*X)\simeq T^*\Sigma$,
and $\Sigma\times_{V}C_V (SS(F))^a$ may be regarded as a subset of
$T^*\Sigma$.
\begin{corollary}
Let $\Sigma$ be a bicharacteristic leaf of $V$. 
Then one has
\eqn
&&SS(\rhom[\shd_X](\shm,{\cal T}
\mu hom(F,{\cal O}_X))\vert_\Sigma)\subset \Sigma\times_{V}C_V (SS(F))^a.
\eneqn
\end{corollary}

\begin{proof}
The map
$\rho_V$ decomposes as
$V\times_{T^*X} T^*(T^*X)\to[j_d] T^*V\to[h] T_V T^*X$.
Here $j:V\hookrightarrow T^*X$ is the embedding.
Set $\shs=\rhom[\shd_X](\shm,{\cal T}
\mu hom(F,{\cal O}_X))$.
Then
the support of $\shs$ is contained in $V$, and
$$SS(\shs\vert_V)=j_d(SS(\shs))\subset h^{-1}(C_V (SS(F))^a).$$
Hence the corollary follows from the following lemma.
\end{proof}

\begin{lemma}
Let $f\colon X\to Y$ be a smooth morphism of real analytic manifolds.
Let $T^*(X/Y)$ be the relative cotangent bundle 
and $h\colon T^*X\to T^*(X/Y)$ be the canonical projection.
Let $S$ be a closed conic subset of $T^*(X/Y)$,
and let $F\in D^b(\C_X)$.
Then the following two conditions are equivalent.
\bnum
\item
$SS(F)\subset h^{-1}(S)$.
\item 
for any $y\in Y$, identifying $f^{-1}(y)\times_X T^*(X/Y)$
with $T^*(f^{-1}(y))$,
$$SS(F\vert_{f^{-1}(y)})\subset f^{-1}(y)\times_X S.$$
\enum
\end{lemma}
\begin{proof}
(i)$\Rightarrow$(ii)\quad
Since
$SS(F\otimes \C_{f^{-1}(y)})=h^{-1}(SS(F\vert_{f^{-1}(y)}))$,
it is enough to show that
$$SS(F\otimes \C_{f^{-1}(y)})\subset h^{-1}S.$$
Since
$SS(F\otimes \C_{f^{-1}(y)})\subset SS(F)\plus T^*_{f^{-1}(y)}X$,
we may reduce the assertion to
\eq\label{eq:2}
h^{-1}S\plus T^*_{f^{-1}(y)}X)\subset h^{-1}S.
\eneq
Let $N=X\times_Y T^*Y\subset T^*X$.
Then $T^*(X/Y)=T^*X/N$, and $T^*_{f^{-1}(y)}X\subset N$.
Hence \eqref{eq:2} is a consequence of
$h^{-1}S\plus N=\ol{h^{-1}S+N}$ (Lemma \ref{L:15} (i)).

\vspace{3pt}
\noindent
(ii)$\Rightarrow$(i)\quad
Let us take a coordinate system $x=(x_1,x_2)$ on $X$ such that $f$ is given
by $x\mapsto x_1$.
Assume that $(x_0,\xi_0)\in T^*X\setminus h^{-1}S$.
Set 
$L(z,\delta,\epsilon)=\{x\in X;
\delta\vert x-z\vert\le-\langle x-z,\xi_0\rangle<\epsilon\delta\}$.
It is enough to show that
$$R\Gamma(X;F\otimes\C_{L(x,\delta,\epsilon)})=0$$
for $x$ sufficiently close to $x_0$ and $0<\epsilon,\delta\ll 1$.
For any $y\in Y$,
$R\Gamma(f^{-1}(y);F\otimes\C_{L(x,\delta,\epsilon)}\vert_{f^{-1}(y)})=0$
by the assumption.
Hence
$Rf_*(F\otimes\C_{L(x,\delta,\epsilon)})=0$,
which implies
$R\Gamma(X;F\otimes\C_{L(x,\delta,\epsilon)})=0$.
\end{proof}

\section{Regularity for $\shd$-Modules}\label{section:reg}

The results contained in this section are extracted or adapted from
\cite{D'A-T} and \cite{K2}.

Recall that a good filtration on a  $\shd_X$-module $\shm$ is  
 a sequence of coherent $\sho_X$-submodules 
$\{\shm_{k}\}_{k\in\Z}$ satisfying:
\bnum
\item
$\shd_X(l) \shm_k\subset {\cal M}_{k+l}$ for any $l,\,k\in \Z$,  
\item
$\shm=\bigcup\limits_k\shm_k$,
\item
locally on $X$,  $\shm_{k}={0}$ for $k\ll 0$,
\item
locally on $X$, $\shd_X(l) \shm_k= \shm_{l+k}$ for $k\gg0$ and
any $l\geq0$.
\enum

We shall use the following notations.
For a coherent $\shd_X$-module $\shm$, we denote by $\she_X \shm$ the
coherent $\she_X$-module  
$\she_X\underset{\pi^{-1}\shd_X}{\otimes}\pi^{-1}\shm$.
For a coherent $\she_X$-module $\shm$ and 
a coherent  $\she_X (0)$-submodule $\shn_0$ of $\shm$,
we set 
$$\shn_0(m):=\she_X (m)\shn_0 \mbox{ for $m\in\Z$.}$$   

\begin{lemma}[\cite{K-K}]\label{L:reg1}
 Let $\shm$ be a coherent $\shd_X$-module and let
$\shn_0$ be a coherent $\she_X(0)$-submodule of 
$\she_X\shm\vert_{\dot{T}^*X}$ such
that $\shn_0$ generates $\she_X\shm$ on $\dot{T}^*X$. Set 
\begin{align*}
\shm_k &=\{u\in \shm; (1\otimes u)\vert_{\dot{T}^*X}\in
\shn_0(k)\}\ \text{ for $k\geq 0$,}\\
\shm_k &=0\ \text{ for $k<0$.}
\end{align*}
Then $\shm_k $  defines a good filtration on 
$\shm$.
\end{lemma}

We shall use the notion of regularity along a closed 
analytic subset $V$  of
$\dot{T}^*X$ due to \cite{K-O}. 
Let us denote by  $\shj_{V}$ the subsheaf of $\she_X$ 
of microdifferential operators of order at most $1$ 
whose  symbol of order $1$ vanishes on $V$. 
In particular $\shj_{V}$ contains $\she_X(0)$. 
Then $\she_V$ denotes the sub-sheaf of rings of
$\she_X$ generated by $\shj_V$. More precisely, 
$$
\she_{V}=\cup_{m\geq0}\shj_V^m.
$$
 
\begin{definition}\cite{K-O}\label{Def:1} 
Let $\shm$ be a coherent $\she_X$-module. An $\she_V$-lattice 
in $\shm$ is an  $\she_V$-submodule 
$\shn_0$ of $\shm$ such that $\shn_0$ is $\she_X(0)$-coherent and 
generates $\shm$ over $\she_X$.
\end{definition}

\begin{definition}\cite{K-O}\label{Def:2}
A coherent $\she_X$-module $\shm$ is called regular along $V$ 
if $\shm$ has an $\she_V$-lattice locally on $\dot{T}^*X$.

Note that if a coherent $\she_X$-module $\shm$ is regular along $V$,
its support is contained in $V$.

If $\shm\in D^b_{\rm coh}(\she_X)$, one says that $\shm$ is regular
along $V$ if $H^k(\shm)$ is regular along $V$ for every $k\in\Z$. 
\end{definition}

\begin{definition}\label{Def:3}
Let $V$ be a closed analytic subset of $T^*X$.
A coherent
$\shd_X$-module $\shm$ is called regular along $V$ 
if the characteristic variety of $\shm$
is contained in $V$ and if 
$\she_X\shm$ is regular along $\dot{V}:=V\setminus T_X^*X$.
\end{definition}
One extends this definition to $D^b_{\rm coh}(\shd_X)$ as in the 
$\she_X$-module case.

The next result follows from \cite{K-O}.
\begin{lemma}\label{L:reg3}
Consider a distinguished triangle $\shm'\to\shm\to\shm''\to[1]$ in 
$D^b_{\rm coh}(\shd_X)$ and assume that two of these three objects are
regular along $V$. Then so is the third one.
\end{lemma}

The next lemma gives a characterization of the
$\shd_X$-modules which are regular along $V$ when $V$ has a special
form.
If $V$  is an involutive vector subbundle of $T^*X$,
then one can find locally on $X$ a smooth morphism 
$f\colon X\to Y$ such that $V=X\times_YT^*Y$
 
\begin{lemma}\label{L:reg2}
Let $f\colon X\to Y$ be  a smooth morphism and let $V:=X\times_YT^*Y$. 
Let $\shm$ be a coherent $\shd_X$-module regular along $V$. 
Then, locally on $X$, 
$\shm$ is the pull-back of a coherent $\shd_Y$-module by $f$.
In particular $\shm$ 
admits a finite resolution  locally on $X$
$$0\to \shl^N\to \shl^{N-1}\to\cdots\to \shl^0\to\shm\to 0,$$  
each $\shd_X$-module $\shl^j$ being
isomorphic to a finite direct sum of the $\shd_X$-module $\shd_{X\to Y}$. 
\end{lemma} 
\begin{proof}
Since the assertion is local on $X$ we may assume that 
$X=Z\times Y$ and  $f$ is the projection $Z\times Y\to Y$.
Hence $\dot{V}=Z\times \dot{T}^*Y$. Let $\Theta_{f} \subset \Theta_{X}$
denote the $\sho_X$-module
of vector fields tangent to the fibers of $f$. With the above
definition, an  $\she_{\dot{V}}$-lattice for
$\she_X\shm\vert_{\dot{T}^*X}$ 
is a coherent $\she_X (0)$-submodule $\shn_0$ of $\she_X\shm$ 
such that $\shn_0$ generates $\she_X\shm$ 
on $\dot{T}^*X$ and such that 
$\Theta_f\shn_0\subset\shn_0$. 

Locally on $X$, there exists a  coherent $\sho_X$-submodule 
$\shm_0$ of $\shm$ such that
$\shm \simeq\shd_X \shm_0$. 
Let us prove that the coherent $\she_{\dot{V}}$-module 
$\shn_0:=\she_{\dot{V}}\shm_0$ is an
 $\she_{\dot{V}}$-lattice in
$\she_X\shm\vert_{\dot{T}^*X}$. 
Since  $\shm_0$ generates $\she_X\shm$ it is
sufficient to prove that  $\shn_0$ is $\she_X(0)$-coherent.

Locally on $X$, there exists a finite covering of $\dot{T}^*X$ by
$\C^\times$-conic open subsets $U_{j}$ and
$\she_{\dot{V}}$-modules $\shn_j$ such that
$\shn_j$ is an $\she_{\dot{V}}$-lattice in $\shm$ on $U_j$.
Hence, for each $j$, there exists $m_{j}\in\Z$ such that 
$\shn_{0}$ is contained in $\shn_j(m_j)$ on $U_j$.

Consider the increasing sequence  
$\{\shj_{V}^{k}\shm_0\}_{k\geq 0}$ 
of coherent $\she_X(0)$-submodules of 
$\shn_0$.
For each $j$, the restriction of this sequence to $U_j$ is contained in
$\shn_j(m_j)$, which is
$\she_X (0)$-coherent, hence it is locally 
stationary. Therefore  $\shn_0$ is $\she_X(0)$-coherent.
 
To summarize, we have constructed an $\she_{\dot{V}}$-module $\shn_0$,
coherent over  $\she_X(0)$-module and which generates $\shm$.
We may now apply Lemma \ref{L:reg1} and consider the  good filtration 
$\shm_{k}$ associated to  $\shn_0$. 

By the construction, $\Theta_f\shm_{k}\subset\shm_k$. 
Setting $\shl=\shm_k$ for $k\gg 0$, we  find a coherent
$\sho_X$-submodule
$\shl$ of $\shm$ such that $\Theta_f\shl
\subset\shl$ and $\shm=\shd_{X}\shl$. Hence,
locally on $X$, there is a $\shd_X$-linear epimorphism 
$\phi_0\colon (\shd_{X\to Y})^{N_0}\to\shm$. 
Repeating the same construction with $\shm$ replaced by $\ker(\phi_0)$, we 
construct an exact sequence 
$$ \shd_{X\to Y}^{N_1}\to[\phi_1]\shd_{X\to Y}^{N_0}\to[\phi_0]
\shm\to 0.$$ 
Since 
$$\hom[\shd_X](\shd_{X\to Y},\shd_{X\to Y})\simeq f^{-1}\shd_Y,$$ 
there is a $\shd_Y$-linear morphism 
$\psi\colon \shd_Y^{N_1}\to\shd_Y^{N_0}$ 
such that
$\phi_1=\opb{f}\psi$. 
Set $\shn=\coker \psi$. 
This is a coherent $\shd_Y$-module and  we have an isomorphism 
$$\shm\simeq \shd_{X\to Y}\tens[f^{-1}\shd_Y]f^{-1}\shn.$$  
To conclude, we choose a  finite  resolution  of
$\shn$ by finitely free
$\shd_{Y}$-modules, and tensorize  over $f^{-1}\shd_{Y}$ by  
the flat $f^{-1}\shd_{Y}$-module 
$\shd_{X\to Y}$.

\end{proof}

\section{Review on  normal cones}\label{section:norcone}

We shall now recall some constructions of \cite{K-S1} which will
be useful for the next steps. To start with, we shall assume that $X$
is a real  manifold. Let $S_1$ and $S_2$ be two subsets of $X$. The
normal cone $C(S_1,S_2)$ is a closed conic subset of
$TX$ which can be described as follows: 

Let $(x)$ be a system of local
coordinates on $X$. Then $(x_0;v_0)\in C(S_1,S_2)$ if and only if there exists a
sequence $\{(x_n,y_n,c_n)\}$ in $S_1\times S_2\times\mathbb{R}^+$ such that 
\begin{equation}\label{E:103}
x_n\underset{n}{\to}x_0,\,y_n\underset{n}{\to}x_0,\,c_n(x_n-y_n)\underset{n}{\to}v_0 
\end{equation}
If
 $A$ is a conic  subset of $T^*X$, we denote by 
$A^{\circ}$ its polar which is  a conic subset  in $TX$. Let $A$ and $B$ be two
conic subsets of
$T^*X$. One  defines the sum 
$$A+B= \{(x;\xi)\in T^*X; \xi= \xi_1+\xi_2\,\, \text{for some}\,\, (x;\xi_1)\in A
\,\,\text{and}\,\,
(x;\xi_2)\in B\}.$$
   If 
$A$ and
$B$ are two closed conic subsets of $T^*X$, one also defines  $A\widehat{+}B$, a closed conic set
containing $A+B$, which may be described in a local canonical coordinate system $(x;\xi)$ as
follows: 
$(x_0;{\xi}_0)$ belongs to $A\widehat{+}B$
if and only if there exist a sequence
$\{(x_n;{\xi}_n)\}_n$ in $A$, a sequence
$\{(y_n;\eta_n)\}_n$ in $B$, such that
\begin{equation}\label{E:10}
x_n\underset{n}{\to}x_0, 
y_n\underset{n}{\to}x_0, 
~(\xi_n+\eta_n)\underset{n}{\to}\xi_0,
|x_n-y_n||\xi_n|\underset{n}{\to}0.  
\end{equation}
 If $A\cap B^{a}\subset T^*_X X$, we have $$ A\widehat{+}B=A+B.$$
Let now $Y$ be a closed submanifold of $X$ and let $i:Y\to X$ be the
 inclusion morphism. Using
 the Hamiltonian isomorphism, we get an embedding of $T^*Y$ 
into $T_{{T}_{Y}^*X}(T^*X)$. Let
$A$ be a conic subset of
$T^*X$. One sets 
$$i^{\sharp}(A):=T^*Y\cap{C_{{T}_{Y}^*X}(A)}.$$ 
This set can be described explicitly by local coordinate systems as
follows.
Let $(x,y)$ be a  local   coordinate system on
$X$ such that  $Y=\{x=0\}$. Then $(y_0;\eta_0)\in
i^{\sharp}(A)\Leftrightarrow
\mbox{there exists a sequence}\,
\{(x_n,y_n;\xi_n,\eta_n)\}_n$ in $A$ such that 
\begin{equation}\label{E:101}
x_n\underset{n}{\to}0,
y_n\underset{n}{\to}y_0,
\eta_n\underset{n}{\to}\eta_0\,\, \mbox{and}\,\,|x_n||\xi_n|\underset{n}{\to}0.
\end{equation}
 
 \begin{lemma}\label{L:15}
 Let $X$ be an open subset in a  finite dimensional real vector space
$E$ with $0\in E$  . Let 
$\Lambda$  be 
     closed  conic  subset of $T^*X$. Let $L$ be a vector subspace of
the dual vector space  $E^*$  and let
$V
=X\times L$. Then we have
\bnum
\item
$V\widehat{+}\Lambda=\overline{V+\Lambda}$.
\item
 For any $\theta\in E^*$ we have
$$
\overline{V+(\Lambda\widehat{+}X\times \mathbb{R}^{\leq 0}\theta)}\cap
X\times \mathbb{R}^{\geq 0}\theta\subset\overline{ V+\Lambda}.$$
\item
Let $i\colon E^*\to T^*X$ be the map $\theta\mapsto(0;\theta)$ and
assume
$(0;0)\in\Lambda$. Then, for any vector subspace $N$ of $E^*$,
$$(\{0\}\times N)\cap(V\widehat{+}\Lambda)\subset T_X^*X$$ if and
only if
$$ (N+i^{-1}(\Lambda))\cap L=\{0\}.$$ 
\enum
\end{lemma} 
\begin{proof}

(i) The inclusion $V\widehat{+}\Lambda\subset\overline{V+\Lambda}$
is clear. Conversely, assume that there are sequences 
 $\{(x_n; \zeta_n)\}_n\subset V$ and 
$\{(y_n;\eta_n)\}_n\subset\Lambda$  
such that $$x_n\underset{n}{\to} x,~ y_n\underset{n}{\to}
x, ~\zeta_n+\eta_n\underset{n}{\to}
\xi.$$ Then the sequence $\{(y_n;\zeta_n)\}_n$ is contained in $V$. Therefore, 
$\{(y_n;\zeta_n)\}_n$ and
$\{(y_n;\eta_n)\}_n$ satisfy (2).

\noindent
(ii) Suppose that $(z;\theta)\in
\overline{V+(\Lambda\widehat{+}X\times\mathbb{R}^{\leq 0}\theta)}$.
 Then by $(i)$ 
there  exist sequences 
$(\zeta_n)_n$ in $L$,  
  $\{(y_n;\eta_n)\}_n$ in $\Lambda$, 
$\lambda_{n} \geq 0,$    such that
$y_n\underset{n}{\to}z$ and
$(\zeta_n +\eta_n - 
\lambda_{n} \theta)\underset{n}{\to}\theta.$ Therefore $$\zeta_n +\eta_n - (
\lambda_{n} +1)\theta\underset{n}{\to} 0.$$ Since $\lambda_n\geq 0$, we get
$$\frac{\zeta_n+\eta_n}{\lambda_n+1}\underset{n}{\to}\theta,$$ hence 
$(z;\theta)\in \overline {V+\Lambda}$.

\noindent
(iii) The condition is obviously necessary. Let us now assume
that there exists $\theta\in N$, $\theta\ne 0$, such that
$(0;\theta)\in V\widehat{+}\Lambda$. Then there exist sequences
$(\zeta_n)_n$ in $L$ and
$(y_n;\eta_n)_n$ in $\Lambda$ such that  $y_n\underset{n}{\to}0$ and $\zeta_n
+\eta_n \underset{n}{\to}\theta$. Taking suitable subsequences, we may 
 assume that  $\zeta_n/|\zeta_n|$ converges to $l\in L, l\ne 0$.
Suppose that
$\zeta_n\underset{n}{\to}0$. Then $(0;\theta)\in\Lambda$, a contradiction.
 If $|\zeta_n|$ is unbounded  we get 
$\zeta_n/|\zeta_n|+\eta_n/|\zeta_n|\underset{n}{\to}0$ hence $(0;-l)\in
L\cap i^{-1}(\Lambda)$, a contradiction.  In the other case,
  we may assume that $\zeta_n\underset{n}{\to}l$ and setting  
$(0;\eta)=(0;\theta-l)\in (N +L)\cap
i^{-1}(\Lambda)$,  we get $\eta+l=\theta$, a contradiction.
 \end{proof}

%

\section{Proof of Theorem \ref{T:1}}
Let $X$ be an open subset of a finite-dimensional real vector
space $E$. 
For $(x_0;\xi_0)\in T^*X$,  
$\epsilon>0$, $\delta>0$, an open convex proper cone
$\gamma$  of $E$ and $v\in\gamma$,
we introduce the following notation:
\eqn
L_{\epsilon,\xi_0}&=&
\{y\in E; \langle y-x_0,\xi_0\rangle >-\epsilon\},\\
Z(x,\gamma,\epsilon,\xi_0)&=&(x+\overline{\gamma})\cap L_{\epsilon,\xi_0}.
\eneqn
Here $\overline{\gamma}$ denotes the closure of $\gamma$, 
The following result is proved in \cite{K-S1}.
\begin{lemma}\label{L:10}
Let $F\in D^b(\C_X)$  and let $p=(x_0;\xi_0)\in T^*X$. 
The conditions below are equivalent:
\bnum
\item
$p\notin SS(F)$.
\item
There exists an open neighborhood $U$ of $x_0$ and an open 
convex proper subanalytic  cone 
$\gamma\subset E$ such that
$\xi_0\in \Int(\gamma^\circ)^a$, satisfying:\\ 
for any $ x\in U$ and sufficiently small $\epsilon>0$,
$Z(x,\gamma,\epsilon,\xi_0)$ is contained in $X$ and
\begin{equation}\label{E:3}  
 \rsect_c (Z(x,\gamma,\epsilon,\xi_0); F)=0.
\end{equation}
\enum      
\end{lemma}

We can now embark into the proof of Theorem \ref{T:1}

We may assume $X=Z\times Y$, $f$ is the projection, and $X,Z,Y$ are
open subsets of affine complex spaces. Moreover, using the results of
\S~\ref{section:reg}, we may assume that $\shm=\shd_{X\to Y}$.
We shall set for short:
$$\shh(F):=\rhom[\shd_X](\shm,{\cal T}hom(F,\sho_X)).$$

Let  $\theta= (x_0;\xi_0)\notin V\plus (SS(F))^a$. 
Let us take an open convex proper cone $\gamma$ 
and an open neighborhood $U$ of $x_0$ such that
$\xi_0\in\text {Int}(\gamma^\circ)^a$ 
and
$U\times \gamma^\circ~\cap \big( V\plus (SS(F))^a\big)\subset T^*_X X$. By 
Lemma \ref{L:10} and keeping its notations, it is enough to prove 
$$ 
\rsect_c\big(Z(x,\gamma,\epsilon,\xi_0);\shh(F)\big)=0
\ \mbox{ for any $x\in U$ and sufficiently small $\epsilon>0$.}
$$
Taking $v\in\gamma$, set
$Z_\delta=Z(x-\delta v,\gamma,\epsilon-\delta,\xi_0)$
for $0<\delta\ll\epsilon$.
Then we have
\[\C_{Z(x,\gamma,\epsilon,\xi_0)}=\indlim[\delta>0]\C_{Z_\delta},\]
and hence we obtain
\begin{eqnarray*}
H^j_c(Z(x,\gamma,\epsilon,\xi_0);\shh(F))&\simeq&
\indlim[\delta>0]H^j_c\bigr(X;
\shh(F)\tens\C_{Z_\delta}\bigr).
\eneqn
Set
\[
Z'_\delta=
(x-\delta v+\gamma)\cap
\{y\in X; \langle y-x_0,\xi_0\rangle \geq-\epsilon+\delta\}.\]
Then for $0<\delta'<\delta$, there is a chain of morphisms
\eqn
\shh(F)\tens\C_{Z_\delta}
\to \shh(F\tens\C_{Z'_\delta})
\to\shh(F)\tens\C_{Z_{\delta'}}
\to \shh(F\tens\C_{Z'_{\delta'}}).
\eneqn
Therefore we have
\eqn
\indlim[\delta>0]H^j_c\bigr(X;
\shh(F)\tens\C_{Z_\delta}\bigr)
&\simeq&
\indlim[\delta>0]H^j_c\bigl(X;\shh(F\tens\C_{Z'_\delta})\bigr).
\end{eqnarray*}
Since $f$ is proper over the support of 
$\C_{Z'_\delta}$, we may apply 
Theorem 7.2 of \cite{K-S2} and obtain  
\begin{eqnarray*}\begin{array}{c}
H^j_c\big(Z(x,\gamma,\epsilon,\xi_0);\shh(F)\big)\simeq
\indlim[\delta>0]H^j_c\Bigl(Y;
{\cal T}hom\bigl(\reim{f}(F\tens\C_{Z'_\delta}),\sho_Y\bigr)\Bigr).
\end{array}
\end{eqnarray*}
Hence, we are reduced to prove
\eq\label{eq:van}
&&\reim{f}(F\tens\C_{Z'(x,\gamma,\epsilon,\xi_0)})=0,
\eneq
where $
Z'(x,\gamma,\epsilon,\xi_0)=
(x+\gamma)\cap
\{y\in X; \langle y-x_0,\xi_0\rangle \geq-\epsilon\}$.
In order to prove this, we shall apply
\cite[Proposition 5.4.17]{K-S1}.
Set $X_t=\{y\in X;\langle y-x_0,\xi_0\rangle\ge t\}$.
Then
if we prove
\eq\label{eq:misup}
&&(y;-\xi_0)\notin \Bigl(SS(F\otimes\C_{(x+\gamma)})+V\Bigl),
\eneq
for $y\in U$,
we have
$$\reim{f}(F\tens\C_{Z'(x,\gamma,\epsilon,\xi_0)})
=\reim{f}(F\tens\C_{(x+\gamma)\cap X_t}).$$
Hence taking $t>0$, we obtain the desired result \eqref{eq:van}.

Thus the proof is reduced to
\eqref{eq:misup}.
We have
\[SS(\C_{(x+\gamma)})\subset
X\times\gamma^\circ{}^a.\]
Hence we have
\[SS(F\otimes\C_{(x+\gamma)})\subset
SS(F)\widehat{+}\bigr(X\times \gamma^\circ{}^a\bigr).\]
Since $SS(F)\cap (X\times \gamma^{\circ})\subset T_X^*X$, we get $$
SS(F)+ (X\times \gamma^{\circ})= SS(F)\widehat{+}
(X\times \gamma^{\circ})$$ and  \[SS(F\otimes\C_{(x+\gamma)})\subset
SS(F)+\bigr(X\times \gamma^\circ{}^a\bigr).\]

On the other hand, by the choice of $\gamma$, we have
\[(X\times \gamma^\circ)\cap
(SS(F)+V)\subset T^*_XX.\]
Hence
\[(X\times \gamma^\circ)\cap
\Bigl(SS(F)+(X\times\gamma^\circ{}^a)+V\Bigr)\subset T^*_XX\]
and we obtain
\[(X\times \Int(\gamma^\circ))\cap
\Bigl(\left(SS(F)\widehat{+}(X\times\gamma^\circ{}^a)\right)+V\Bigr)
=\emptyset.\]
Then the desired result follows from 
 $-\xi_0\in \Int(\gamma^\circ)$.

\section{Proof of Theorem \ref{T:2}}


We shall now embark in the proof of Theorem \ref{T:2}.
 Since the question is local on $X$, by Lemma 3.6  we may assume that
${\cal M}$ is isomorphic to
${\cal D}_{X\to Y}$. By (7.5) of \cite{K-S2}, if $d$ denotes the codimension of $Z$, we have a natural
isomorphism
\begin{equation}\label{E:36}
R{\cal H}om_{{\cal
D}_Z}({\cal M}_Z, {\cal T} hom(F_Z, {\cal O}_Z))\simeq{R{\cal
H}om_{{\cal D}_X}({\cal M}, {\cal T }hom(F_Z, {\cal
O}_X))|_Z}[2d]\end{equation} With this isomorphism in hand, Theorem
2.2 will be  a consequence of the next Lemma for a real submanifold.

If $Z$ is  a real submanifold of a complex manifold $X$, we still denote by $T^*_Z X$ the
conormal bundle to $Z$ of the real underlying manifold $X^{\mathbb{R}}$ . 
\begin{lemma}\label{L:35}
Let $Z$ be a real analytic closed  submanifold of $X$ of codimension $d\geq1$ and assume that
$$T^*_Z X\cap( V\widehat{+} SS(F)^a) \subset T^*_X X.$$ Let $\cal{M}$ be a
 regular ${\cal
{D}}_X$-module along $V$. Then  the following natural morphism is an
isomorphism:\begin{equation}\label{E:37}
R{\cal H}om_{{\cal D}_X}({\cal M}, {\cal T}hom(F,{\cal
O}_X))|_Z\to{R{\cal H}om_{{\cal D}_X}({\cal M}, {\cal T}hom(F_Z,
{\cal O}_X))|_Z}[d]\end{equation}
\end{lemma}
\begin{proof}
Since $T_Z^*X\cap V\subset T_X^*X$,  there exists a
coordinate system $(y_1,\cdots,y_N)$ on  $Y$ as a real analytic
manifold and a coordinate system  $
(f_1,\cdots,f_m,y_1\circ f,\cdots,y_N\circ f)$ of $X$ such that  
$Z$ is defined by the equations
$f_1 =\dots =f_d =0$.
 We shall argue by induction on $d$.\\$(i)$ Let us prove the result for  $d=1$. 
 Assume that $Z$ is a real analytic hypersurface defined by the equation
$f(z)=0$, and set 
$${Z}^{-}=\{z\in X; f(z)\leq0\},\,\,\,   {Z}^{+}=\{z\in
X; f(z)\geq0\}.$$ Assume that $df, -df\notin V\widehat{+}SS(F)^a$. We
shall show that the morphism 
\begin{equation}\label{E:35}  R{\cal H}
om_{{\cal D}_X}({\cal M}, {\cal T}hom(F,{\cal O}_X))|_Z
\to{R{\cal H}om_{{\cal
D}_X}({\cal M}, {\cal T}hom(F_Z, {\cal O}_X))|_Z}[1].
\end{equation}
is an isomorphism. The  morphism (\ref{E:35}) is given by $$ F\to
F_{Z^+}\oplus F_{Z^-}\to F_Z\overset{+1}{\to}$$
  Therefore to obtain (\ref{E:35}) it is enough to prove that $$R{\cal
H} om_{{\cal D}_X}({\cal M},{\cal T}hom (F_{ Z^{\pm}},{\cal
O}_X))|_Z=0.$$  Since $$SS(F_{Z^{\pm}})\subset
SS(F)\widehat{+}(X\times \mathbb{R}^{\geq 0}(\pm df)),$$
$\pm df\in
V\hat{+}SS(F_{ {Z}^{\pm}})^a$ implies 
 $\pm df\in V\hat{+}SS(F)^a$  by Lemma \ref{L:15} (ii),
which contradicts the assumption. Hence         
$$\pm df\notin V\hat{+}SS(F_{ {Z}^{\pm}})^a.$$ 
Denoting  
 $$ {\cal S}_{\pm} ={\cal H}
om_{{\cal D}_X}({\cal M},{\cal T}hom (F_{{Z}^{\pm}},{\cal
O}_X)),$$ then, by Theorem \ref {T:1}, $\pm df\notin SS({\cal
S}_{\pm})$ hence
${\cal S}_{\pm}|_Z=R\Gamma _{Z^{\pm}}({\cal
S}_{\pm}|_Z)=0$.

\noindent $(ii)$ Let us set 
$Z_1=\{f_1=0\}$ and $Z_2=\{f_2
=\dots=f_d=0\}$.  Since $$T^*_{Z_2}X\cap
(V\widehat{+}(SS(F)^a))\subset T^*_X X$$ in a neighborhood of $Z$, the
hypothesis of induction implies 
$$ R{\cal H}om_{{\cal D}_X}({\cal M}, {\cal T}hom(F,{\cal
O}_X))|_Z \simeq {R{\cal H}om_{{\cal D}_X}({\cal M},
{\cal T}hom(F_{Z_{2}}, {\cal O}_X))|_Z}[d-1].$$
On the other hand, $SS(F_{Z_2})\subset
SS(F)\widehat{+}T^*_{Z_2}X$.
 Hence by $(i)$ 
$$ R{\cal H}om_{{\cal D}_X}({\cal M},{\cal T}hom(F,{\cal
O}_X))|_Z\simeq {R{\cal H}om_{{\cal D}_X}({\cal M},{\cal T}hom
(F_{Z_{2}}, {\cal O}_X))|_Z}[d-1]$$
$$\simeq {R{\cal
H}om_{{\cal D}_X}({\cal M}, {\cal T}hom(F_{Z}, {\cal
O}_X))|_Z}[d].$$
 
This ends the proof of (\ref{E:37}).

\end{proof}

\section{Proof of Theorem 2.3 }\label{S:3}
 We start by
recalling the functor of tempered microlocalization.

\subsection{Review on ${\cal T}\mu hom$}\label{SS:1}

We shall recall the construction
 of the functor ${\cal T}\nu hom(\cdot,{\cal O}_X)$ of tempered
specialization of (\cite{A}).

Let $\tilde{X}^{\mathbb{C}}$ be the complex normal deformation of $X\times{X}$ along
 the diagonal $\Delta$ which we
 identify with $X$ by the first
 projection $p_1$.  We may then     
 identify
$TX$ with the normal bundle $T_{\Delta}(X\times{X})$. Let
$t:\tilde{X}^{\mathbb{C}}\to{\mathbb{C}}$ and
$p:\tilde{X}^{\mathbb{C}}\to{X\times{X}}$ be  the
 canonical maps,
let $\tilde{\Omega}$ be $t^{-1}({\mathbb{C}}-\{0\})$ and
$\Omega=t^{-1}(\mathbb{R}^{+})\subset\tilde{\Omega}$.
 Let
$p_2:X\times{X}\to X$ be the  second projection.
 
Consider the following diagram of morphisms:
\begin{equation}\label{E:7}
TX\simeq{T_{\Delta}(X\times{X})}\overset{i}{\hookrightarrow}\tilde{X}^{\mathbb{C}}\overset{j}{\hookleftarrow}\Omega=t^{-1}(\mathbb{R}^{+})  
\end{equation}
Let $\tilde{p}:\tilde{\Omega}\to{X\times{X}}$,~
 be the
 restriction of $p$. Finally denote by $\overline{p}_{1}$ the composition $p_{1}\circ{p}$ and by
$\overline{p}_{2}$ the composition $p_{2}\circ{p}$. 

Under these notations, ${\cal T}\nu hom(F,{\cal O}_X)$ is defined by 
 $${\cal T}\nu hom(F,{\cal O}_X)
={i^{-1}R{\cal H}om_{{\cal D}_{\tilde{X}^{\mathbb{C}}}}({\cal
D}_{\tilde{X}^{\mathbb{C}}\underset{\ol{p}_{1}}{\to}X}
,{\cal T}hom(\ol{p}_{2}^{-1}F\otimes{\mathbb{C}_{\Omega}},{\cal
O}_{\tilde{X}^{\mathbb{C}}}))}.$$ 
Let $D^b_{\R^+}(\mathbb{C}_{TX})$
(resp. $D^b_{\R^+}(\mathbb{C}_{T^*X})$)
be the derived category
of complexes of sheaves on $TX$ (resp. $T^*X$)
with conic cohomologies.
We denote by the symbol $\,\,\widehat{}\,\,$ the
Fourier-Sato Transform from
$D^b_{\R^+} (\mathbb{C}_{TX})$ to
$D^b_{\R^+} (\mathbb{C}_{T^*X})$. Then
by definition, ${\cal T}\mu hom(F,{\cal
O}_X)={\cal T}\nu hom(F,{\cal O}_X)^{\,\widehat{}}$ .
 Let us recall that under the identification of $T^*(TX)$ with
$T^*(T^*X)$ by  the Hamiltonian isomorphism  
we have  $SS(F)=SS(F)^{\,\widehat{}}$ for any
$F\in{D^b_{\R^+}(\mathbb{C}_{TX})}$.

  Remark that for any coherent ${\cal
D}_X$-module
${\cal M}$, one has
 \begin{align}\label{E:49}
 & R{\cal H}om_{\pi^{-1}{\cal D}_X}(\pi^{-1}{\cal
M},{\cal T} \mu hom(F,{\cal O}_X))\\
&\simeq
{ R{\cal
H}om_{\tau^{-1}{\cal D}_X}(\tau^{-1}{\cal M},{\cal T}\nu hom(F,{\cal
O}_X))}^{\,\widehat{}}.\notag
\end{align}

The proof of Theorem 2.3  will be performed in two steps.

\subsection{First reduction}
 First of all remark  that the statement of Theorem \ref{T:40} is local
on $T^*X$  and invariant by local canonical transformation as proved
in \cite{A},  5.5.  
 Therefore, since $V$ is regular involutive, locally in $\dot{T}^*X$,
we may
 choose a canonical coordinate system $(x;\xi)$, $x=(x_1,...,x_n)$, $\xi=(\xi_1,...,\xi_n)$, 
such that $V=\{(x;\xi)\in{\dot{T}^*X};\xi_1=0,...,\xi_d=0\}$, in other
words, locally on X, we have
$X\simeq Z\times Y$ where $Y$ is an open
subset of $\mathbb{C}^{n-d}$ and 
$V=X\underset{Y}{\times}T^*Y $,
 is associated to the projection $f:X\to Y$.    By the
results
 of \cite{K-O},  we have an exact
sequence of coherent
${\cal E}_X$-modules
 \begin{equation}\label{E:53}
{0}\to {\cal N}\to ({\cal E}_{X}
/({\cal E}_{X}D_{x_{1}}+{\cal E}_{X}D_{x_{2}}+...+{\cal E}_{X}D_{x_{d}}))^N \to{\cal M}\to{0}
\end{equation} 
where ${\cal N}$ is still regular along $V$.

 By
  "devissage" thanks to   (\ref{E:53}), we may then assume that

${\cal M}={\cal E}_X /({\cal E}_{X}D_{x_{1}}+{\cal E}_{X}D_{x_{2}}+...+{\cal
E}_{X}D_{x_{d}})
={\cal E}_{X\to Y}$.\\


 Hence from now on we will assume ${\cal M}={\cal
E}_{X\to Y}$. Of course, in that case, ${\cal M}\simeq{{\cal E}_{X}\otimes_{\pi^{-1}{\cal
D}_X}{\pi^{-1}{\cal D}_{X\to Y}}}$.
\\
\subsection{ End of the proof}\label{SS:2}
 Having (\ref{E:49}) in mind,
 we are bound to prove the analogue of (\ref{E:52}) with ${\cal T}\mu
hom(F,{\cal O}_X)$  replaced by ${\cal T}\nu hom(F,{\cal O}_X)$.
We have

\begin{align}&R{\cal H}om_{\tau^{-1}{\cal
D}_X}(\tau^{-1}{\cal M},{\cal T}\nu hom(F,{\cal O}_X))\notag\\
&\simeq
{i^{-1}R{\cal H}om_{{\cal D}_{\tilde{X}}}({\cal D}_{\tilde{X}^{\mathbb{C}}\underset{h}
{\rightarrow}Y},{\cal
T}hom(\ol{p}_{2}^{-1}F\otimes{\mathbb{C}_{\Omega}},{\cal
O}_{\tilde{X}^{\mathbb{C}}}))}\notag
\end{align}
\noindent where $h=f\circ{\overline{p}_1}$ is a smooth morphism since $f$ and $\overline{p}_1$ are
smooth.

 By Proposition 6.6.2 of \cite{K-S1}, we have an 
inclusion
\begin{align}\label{E:55}
&SS(~i^{-1}R{\cal H}om_{{\cal D}_{\tilde{X}}^{\mathbb{C}}}({\cal
D}_{\tilde{X}^{\mathbb{C}}\underset{h}
{\rightarrow}Y},{\cal
T}hom(\ol{p}_{2}^{-1}F\otimes{\mathbb{C}_{\Omega}},{\cal
O}_{\tilde{X}^{\mathbb{C}}})))\\ &\subset{i^{\sharp}(SS(R{\cal
H}om_{{\cal D}_{\tilde{X}^{\mathbb{C}}}}({\cal
D}_{\tilde{X}^{\mathbb{C}}
\underset{h}
{\rightarrow}Y},{\cal
T}hom(\ol{p}_{2}^{-1}F\otimes{\mathbb{C}_{\Omega}},{\cal
O}_{\tilde{X^{\mathbb{C}}}}))))}\notag.
\end{align}
Therefore it is enough to consider the case of the partial De Rham
system ${\cal D}_{\tilde{X}^{\mathbb{C}}\underset{h}
{\rightarrow}Y}$. 

Let us take a local coordinate system $(x)=(x_1,x_2)$ on $X$
such that $f\colon X\to Y$ is by $(x_1,x_2)\mapsto x_1$.
Then $V=X\times_YT^*Y$ is given by
$V=\{(x_1,x_2;\xi_1,\xi_2);\xi_2=0\}$.
Endow $X\times{X}$
with the system of local coordinates 
$(x,x')$, so that
 $\Delta\subset{X\times{X}}$ is defined by $x=x'$. Under 
the change of coordinates : $x$, $y=x-x'$, $\Delta$ will be
defined by $y=0$. Using $(x,y)$, $\tilde{X}^{\mathbb{C}}$ is
endowed with the coordinates $(t,x,y)$,
and
\begin{equation}\label{E:51} 
\mbox{$p(t,x,y)
=(x,x-ty)$,  $\overline{p}_1(t,x,y)=x$, and 
$\overline{p}_2(t,x,y)=x-ty$}
\end{equation}
Let $(t,x,y;\tau,\xi,\eta)$ be the associated coordinates of 
$T^*(\tilde{X}^{\mathbb{C}})$.

Let  $\tilde{V}$  be the submanifold
${\ol{p}_1}_d(\tilde{X}^{\mathbb{C}}\underset{Y}{\times}T^*Y)$
of $T^*\tilde{X}^{\mathbb{C}}$,  which is explicitly given by 
$$\tilde{V}=\{(t,x,y;\tau,\xi,\eta);\tau=0,\eta=0,\,(x;\xi)\in V\}.$$
By Theorem \ref{T:1} we have the following estimate:
\begin{equation}\label{E:56}
SS(R{\cal H}om_{{\cal D}_{\tilde{X}^{\mathbb{C}}}}
({\cal D}_{\tilde{X}^{\mathbb{C}}
\underset{h}
\to{Y}},{\cal T}hom(\ol{p}_{2}^{-1}F
\otimes{\mathbb{C}_{\Omega}},{\cal
O}_{\tilde{X}^{\mathbb{C}}})))
\subset{\tilde{V}\hat{+}SS(~\overline{p}_{2}^{-1}F
\otimes\mathbb{C}_{\Omega}})^a.
\end{equation}
 By  (\ref{E:55}) and (\ref{E:56}) we get
 \begin{equation}\label{E:57}
SS(R{\cal H}om_{\tau^{-1}{\cal D}_X}(\tau^{-1}{\cal M},{\cal T}\nu
hom(F,{\cal O}_X))) \subset{i}^{\sharp}(  
\tilde{V}\widehat{+}SS(\overline{p}_{2}^{-1}F\otimes
\mathbb{C}_{\Omega})^a).
\end{equation}
Therefore it is enough to prove   the inclusion
\[{i}^{\sharp}(\tilde V\widehat{+}SS(~\overline{p}^{-1}_2
F\otimes\mathbb{C}_{\Omega})^a)\subset{C(V,SS( F))}.\]
More precisely,
since $SS(\overline{p}_{2}^{-1}F\otimes\mathbb{C}_{\Omega})\subset
SS(\overline{p}_{2}^{-1}F)\widehat{+}SS(\mathbb{C}_{\Omega})$ we
shall prove the inclusion 
\[{i}^{\sharp}(\tilde V\widehat{+}(SS(~\overline{p}^{-1}_2
F)^a\widehat{+}SS(\mathbb{C}_{\Omega})^a))\subset{C(V,SS( F))}.\]  We
have 
\begin{align}
 SS(\mathbb{C}_{\Omega}) &= \{(t,x,y;\tau,\xi,\eta);
\xi=0, \eta=0,
\rm{Im}\, t=0,\rm{Re}\, t\geq0,\rm{Re}\,\tau=0\}\\
&\cup
\{(t,x,y;\tau,\xi,\eta);\xi=0, \eta=0,
t=0,\rm{Re}\,\tau\leq 0\},\notag
\end{align}
 and since
$\overline{p}_2$ is smooth,
\begin{align} 
SS(\overline{p}_{2}^{-1}F)&
=\overline{p}_{2,d}\overline{p}_{2,\pi}^{-1}(SS(F))\\
&=\{(t,x,y;\tau,\xi,\eta); (x-ty;\xi)\in SS(F),
 &\eta=-t \xi,\, \tau=-\langle \xi,y\rangle\}.\notag
\end{align}
Hence $SS(\overline{p}_{2}^{-1}F)^a\cap  SS(\mathbb{C}_{\Omega})\subset
T_X^*X$, which implies
$$SS(~\overline{p}^{-1}_2F)^a\widehat{+}SS(\mathbb{C}_{\Omega})^a
=SS(~\overline{p}^{-1}_2F)^a{+}SS(\mathbb{C}_{\Omega})^a.$$
Remark that the identification $T^*(TX)$ with $T(T^*X)$ is
described by 
\[T^*(TX)\ni(x,y;\xi,\eta)\longleftrightarrow (x,\eta;y,\xi)\in T(T^*X).\]
Let $(x_0,y_0;\xi_0,\eta_0)\in T^*(TX)$ and assume that 
$$
(x_0,y_0;\xi_0,\eta_0)\in {i}^{\sharp}(\tilde V
\widehat{+}(SS(~\overline{p}^{-1}_2
F)^a{+}SS(\mathbb{C}_{\Omega})^a)).$$
Then there exist sequences

\noindent $\{(t_n,x_n,y_n;0,\xi_n,0)\}_{n}$ in $ \tilde{V}$,
$\{(t_n',x_n',y_n';\tau_n',\xi_n',\eta_n')\}_n$ in
$SS(\overline{p}_{2}^{-1}F)^a$

\noindent and
$\{(t_n'',x_n'',y_n'';\tau_n'',0,0)\}_n$ in
$SS(\mathbb{C}_{\Omega})^a$ such that
\bnum
\item
$t_n\underset{n}{\to}0$,
$t_n'\underset{n}{\to}0$,
$t_n''\underset{n}{\to}0$,
$t'_n, t''_n\geq0$.
\item
$x_n\underset{n}{\to} x_0$, $x_n'\underset{n}{\to} x_0$,
$x_n''\underset{n}{\to} x_0$.
\item $y_n\underset{n}{\to}y_0$,
$y_n'\underset{n}{\to} y_0$, $y_n''\underset{n}{\to} y_0$.
\item
$\tau_n'+\tau_n''\underset{n}{\to} 0$.
\item
$\xi_n+\xi_n'\underset{n}{\to}\xi_0$.
\item
$\eta_n'\underset{n}{\to}\eta_0$ (hence
$t'_n\xi'_n\underset{n}{\to}-\eta_0$, and 
$t'_n\xi_n\underset{n}{\to}\eta_0$ by (v)).
\enum
By (vi), there exists
a sequence of
positive numbers $(a_n)_n$ 
such that $a_n\underset{n}{\to}0$
and $a_n\xi_n\underset{n}{\to} \eta_0$.
Consider the sequence
$\{(x'_n-t'_ny'_n;-a_n\xi'_n)\}_n$ in $SS(F)$ and
$\{(x'_n+(a_n-t'_n)y'_n;a_n\xi_n)\}_n$ in $ V$. 
Then
$$\mbox{$(x'_n-t'_ny'_n;-a_n\xi'_n)\underset{n}{\to}(x_0;\eta_0)$,
$\{(x'_n+(a_n-t'_n)y'_n;a_n\xi_n)\underset{n}{\to}(x_0;\eta_0)$}$$
and
$${a_n}^{-1}\bigl((x'_n+(a_n-t'_n)y'_n,a_n\xi_n)
-(x'_n-t'_ny'_n,-a_n\xi'_n)\bigr)
=(y'_n,\xi_n+\xi'_n)\underset{n}{\to}(y_0,\xi_0).$$
Hence one has $(x_0,\eta_0;y_0,\xi_0)\in C(V, SS(F))$.

Since $C(V,SS(F))={\rho_V}^{-1}(C_V(SS(F))^a)$,
we finally obtain
 \[ SS(~R{\cal H}om_{\pi^{-1}{\cal D}_X}(\pi^{-1}{\cal
M},{\cal T}\mu hom(F,{\cal O}_X)))
\subset {\rho_V}^{-1}(C_V(SS(F))^a)
\] as asserted.
\qed

\newpage
\small
Teresa Monteiro Fernandes\\
Centro de \'Algebra da Universidade de Lisboa,
Complexo 2,\\
2 Avenida Prof.Gama Pinto, 1699 Lisboa codex Portugal\\
tmf@ptmat.lmc.fc.ul.pt\\
\vspace{15mm}

Masaki Kashiwara\\
Research Institute for Mathematical Sciences,\\
Kyoto University\\
Kyoto 606-8502\\
Japan\\
masaki@kurims.kyoto-u.ac.jp\\
\vspace{15mm}

Pierre Schapira\\
Universit\'e Pierre et Marie Curie, case 82\\
Analyse Alg\'ebrique, UMR7586\\
4, place Jussieu,75252 Paris cedex 05\\
France\\
schapira@math.jussieu.fr\\
http://www.math.jussieu.fr/{\~{}} schapira/\\

\end{document}